\documentclass{article}
\usepackage{amsfonts}
\usepackage{amssymb}
\usepackage{amsmath}
\usepackage{ragged2e}

\begin{document}

\bigskip

\bigskip

\begin{center}
\textbf{THE LIE - TROTTER INTEGRATOR IN THE DYNAMICS OF THE SYMMETRIC FREE RIGID
BODY }
\end{center}
\begin{center}
\bigskip
\end{center}
\begin{center}
\textbf{Ciprian Hedrea}

\bigskip

\textit{Department of Mathematics}

\textit{Polytechnic University of Timi\c{s}oara}

\textit{Square of Victory no. 2, 300006}

\textit{Timi\c{s}oara, Rom\^{a}nia}

\textit{ciprian.hedrea@upt.ro}

\end{center}

\bigskip

\bigskip

\bigskip

\textbf{Abstract.} The numerical integration plays a fundamental role in understanding the behaviour of many mechanical systems. In this paper some important aspects of the mechanical integrators on the dynamics of a mechanical system are studied. More specific, we have shown that if that the Lie-Trotter integrator is obtained, in case of Euler equations for the dynamics of symmetric free rigid body, then it is a Poisson integrator. At the end of the paper some important remarks are presented. 
\bigskip

\textit{Keywords}: Integrator; rigid body; symmetrical rigid body; mechanical system; dynamics.

\bigskip
2000 Mathematics Subject Classification: 37M15, 70Q05 

\bigskip

\textbf{1. Introduction:}

\ \ \

The overall situation for mechanical integrators is a complex one and it is still evolving. 
There are numerical integration methods that preserve some of the invariants of the mechanical system, such as: energy, momentum, or the symplectic form. 

It is well known that if the energy and momentum map include all the integrals from a certain class that one cannot create integrators that are symplectic, energy preserving and momentum preserving unless they coincidentally integrate the equations exactly up to a time parametrization [6].
In mechanics the numerical integration plays a fundamental role.

Because
there are many mechanical systems whose explicit integration is unknown,
in our effort to understand their behavior, it is very useful the
numerical integration of their dynamics.

Recently, more research is dedicated to the use of the mechanical integrators in molecular dynamics, spin systems, magnetism, etc.

\bigskip

\textbf{2. Numerical and symplectic integrators}

\ \ \ 

\textbf{Definition 2.1.} A numerical integrator on $%
\mathbb{R}
^{n}$ consists in a number of applications differentiable on class $C^{\infty }$%
\begin{equation*}
\phi _{t}:%
\mathbb{R}
^{n}\rightarrow 
\mathbb{R}
^{n},
\end{equation*}%
\ which are differentiable in relation to $t\in 
\mathbb{R}
.$ We write 
\begin{equation*}
x_{n+1}=\phi _{t}\left( x_{n}\right) .
\end{equation*}

\bigskip

We consider in $%
\mathbb{R}
^{n}$ the system of differential equations given by%
\begin{equation}
\overset{\cdot }{x}\ =f\left( x\right) ,  \tag{2.1}
\end{equation}%
where $x\in 
\mathbb{R}
^{n}$ and $f\in C^{\infty }\left( 
\mathbb{R}
^{n},%
\mathbb{R}
^{n}\right) .$

In the following, we will present a series of numerical integrators that enable us to
approximate the solutions of the system $\left( 2.1\right) $.

$\left( i\right) $ The Gauss-Legendre integrator, given by%
\begin{equation}
\frac{x_{n+1}-x_{n}}{h}=f\left( \frac{x_{n+1}+x_{n}}{2}\right) .  \tag{2.2}
\end{equation}

$\left( ii\right) \ $The Runge-Kutta integrator with $s$ steps, given by%
\begin{equation}
\left\{ 
\begin{array}{c}
\frac{x_{n+1}-x_{n}}{h}=\sum\limits_{i=1}^{s}b_{i}\ f\left( X_{i}\right) ,
\\ 
X_{i}=x_{n}+\sum\limits_{j=1}^{s}a_{ij}\ f\left( X_{j}\right) .%
\end{array}%
\right.  \tag{2.3}
\end{equation}

$\left( iii\right) \ $The Euler integrator, given by%
\begin{equation}
\frac{x_{n+1}-x_{n}}{h}=f\left( x_{n}\right) .  \tag{2.4}
\end{equation}

$\left( iv\right) \ $The modified Euler integrator, given by%
\begin{equation}
\frac{x_{n+1}-x_{n}}{h}=f\left( x_{n}\right) +f\left( x_{n+1}\right) . 
\tag{2.5}
\end{equation}

(For other examples, see [3]  and [9]).

\bigskip

Let $\left( M,\omega ,H\right) $ be a mechanical Hamiltonian system.

\textbf{Definition 2.2.} A numerical integrator $\left\{ \phi _{t}\right\}
_{t\in 
\mathbb{R}
}$ on $M$ is called symplectic integrator if%
\begin{equation*}
\phi _{t}^{\ast }\omega =\omega ,\ \left( \forall \right) t\in 
\mathbb{R}
.
\end{equation*}

\textbf{Example 2.1.} Considering a mechanical Hamiltonian system $\left(
M,\omega ,H\right) $, where we have%
\begin{equation*}
M=T^{\ast }%
\mathbb{R}
\simeq 
\mathbb{R}
^{2},\ \omega =dp\wedge dq,\ H\left( p,q\right) =\frac{1}{2}\left(
q^{2}+p^{2}\right) .
\end{equation*}

Then the Ruth integrator associated with this dynamic%
\begin{equation*}
\left\{ 
\begin{array}{c}
\overset{\cdot }{q}=p, \\ 
\overset{\cdot }{p}=-q,%
\end{array}%
\right.
\end{equation*}%
will be given as%
\begin{equation}
\left\{ 
\begin{array}{c}
x_{n+1}^{1}=x_{n}^{1}+hx_{n}^{2}, \\ 
x_{n+1}^{2}=-hx_{n}^{1}+\left( 1-h^{2}\right) x_{n}^{2}.%
\end{array}%
\right.  \tag{2.6}
\end{equation}

The conclusions that can be drawn immediately, are

$\left( i\right) \ $The Integrator $\left( 2.6\right) $ doesn't conserve
energy, because

\begin{equation*}
H\left( x_{n+1}^{1},x_{n+1}^{2}\right) =\left( 1+h^{2}\right) H\left(
x_{n}^{1},x_{n}^{2}\right) +h^{2}\left( \frac{h^{2}}{2}-1\right)
x_{n}^{2}+h^{3}x_{n}^{1}x_{n}^{2}.
\end{equation*}

$\left( ii\right) $ The integrator $(2.6)$ is a symplectic integrator

More exactly we have%
\begin{equation*}
dx_{n+1}^{2}\wedge dx_{n+1}^{1}=dx_{n}^{2}\wedge dx_{n}^{1}.
\end{equation*}

\textbf{Proposition 2.1. }([5]) Let the Hamiltonian mechanical system 
\begin{equation*}
\left( 
\mathbb{R}
^{2n},\omega ,H\right) ,
\end{equation*}%
where%
\begin{equation*}
\omega =dp_{1}\wedge dq^{1}+...+dp_{n}\wedge dq^{n}.
\end{equation*}

Then the corresponding Gauss - Legendre integrator is a symplectic integrator.

\textbf{Proposition 2.2. }([7], [15], [17]) Let the Hamiltonian
mechanical system $\left( 
\mathbb{R}
^{2n},\omega ,H\right) ,$ where%
\begin{equation*}
\omega =dp_{1}\wedge dq^{1}+...+dp_{n}\wedge dq^{n}.
\end{equation*}

Then the Runge - Kutta integrator with $s$ stairs and $h$ step is a
symplectic integrator if and only if%
\begin{equation}
\left\{ 
\begin{array}{c}
b_{i}a_{ij}+b_{j}a_{ji}-b_{i}b_{j}=0, \\ 
1\leq i,\ j\leq s.\ \ \ \ \ \ \ \ \ \ \ \ \ \ \ \ \ 
\end{array}%
\right.  \tag{2.7}
\end{equation}

\textbf{Definition 2.3. }Let $\left( M,\omega ,H,G\right) $ be a mechanical
Hamiltonian system with symmetry having an application moment%
\begin{equation*}
J:M\rightarrow \mathcal{G}^{\ast }.
\end{equation*}

An integrator $\left\{ \phi _{t}\right\} _{t\in 
\mathbb{R}
}$ on $M$ is called integrator moment if%
\begin{equation*}
\phi _{t}^{\ast }J=J,\ \left( \forall \right) \ t\in 
\mathbb{R}
.
\end{equation*}

\bigskip

\textbf{Proposition 2.3. }([6]) Let $\left( M,\omega ,H,G\right) $ be a
mechanical Hamiltonian system with symmetry and having a moment application%
\begin{equation*}
J:M\rightarrow \mathcal{G}^{\ast }.
\end{equation*}

Then any integrator $\left\{ \phi _{t}\right\} _{t\in 
\mathbb{R}
}\ $on $M$ wich is $G$-invariant and is symplectic, is a moment integrator.

\bigskip

\textbf{Proposition 2.4. }([6]) Let $\left( M,\omega ,H,G\right) $ a
mechanical Hamiltonian system with symmetry and having a moment application%
\begin{equation*}
J:M\rightarrow \mathcal{G}.
\end{equation*}

Then if $\left\{ \phi _{t}\right\} _{t\in 
\mathbb{R}
}$ is an integrator on $M$ such that 
\begin{eqnarray*}
&&\left( i\right) \ \left\{ \phi _{t}\right\} _{t\in 
\mathbb{R}
}\text{ is an symplectic integrator,} \\
&&\left( ii\right) \ \left\{ \phi _{t}\right\} _{t\in 
\mathbb{R}
}\text{ is an energy integrator,} \\
&&\left( iii\right) \ \left\{ \phi _{t}\right\} _{t\in 
\mathbb{R}
}\text{ is an moment integrator,} \\
&&\left( iv\right) \ \text{dynamics on the reduced space }\left(
M_{0},\omega _{0}\right) \text{ not integrable.}
\end{eqnarray*}

Then $\left\{ \phi _{t}\right\} _{t\in 
\mathbb{R}
}$ gives the exact solution for initial system dynamics (unreduced).

\bigskip

\textbf{Observation 2.1.} Higher order integrators have generally similar
properties. For details see [14], [10], [4], [19] and [16].

\bigskip

\bigskip

\textbf{3. Poisson integrators}

\ \ \ 

Let $\left( P,\left\{ \cdot ,\cdot \right\} \right) $ be a dimensional finite Poisson manifold
from class $C^{\infty }.$

\textbf{Definition 3.1.} An integrator $\left\{ \phi _{t}\right\} _{t\in 
\mathbb{R}
}$ on $P$ is called Poisson integrator if%
\begin{equation}
\phi _{t}^{\ast }\left\{ f,g\right\} =\left\{ \phi _{t}^{\ast }f,\phi
_{t}^{\ast }g\right\} ,\ \ \left( \forall \right) \ f,g\in C^{\infty }\left(
P,%
\mathbb{R}
\right) .  \tag{3.1}
\end{equation}

\textbf{Observation 3.1. }In the particular case $P=%
\mathbb{R}
^{n},$ relation $(8)$ is equivalent to%
\begin{equation}
\left\{ 
\begin{array}{c}
\left( D\phi _{t}\left( x\right) \right) \cdot \Pi \cdot \left( D\phi
_{t}\left( x\right) \right) ^{t}=\Pi \left( \phi _{t}\left( x\right) \right)
, \\ 
\left( \forall \right) \ t\in 
\mathbb{R}
,\ \left( \forall \right) \ x\in 
\mathbb{R}
^{n},\ \ \ \ \ \ \ \ \ \ \ \ \ \ \ \ \ \ \ \ \ 
\end{array}%
\right.  \tag{3.2}
\end{equation}%
where $D\phi _{t}$ designate Fr\'{e}ch\'{e}t differential of $\phi _{t}$,
and $\Pi $ is the matrix associated to the Poisson structure $\left\{ \cdot
,\cdot \right\} ,$ hence $\Pi =\left[ \left\{ x_{i},x_{j}\right\} \right] .$

\textbf{Example 3.1. }([11]) Let $\left( 
\mathbb{R}
^{2},\Pi ,H\right) $ be a Hamilton-Poisson mechanical system , where%
\begin{equation*}
\Pi =\left[ 
\begin{array}{cc}
0 & x_{2} \\ 
-x_{2} & 0%
\end{array}%
\right] ,
\end{equation*}%
and%
\begin{equation*}
H\left( x_{1},x_{2}\right) =Ax_{1}+Bx_{2}+C,
\end{equation*}%
$A,B,C\in 
\mathbb{R}
.$ A simple calculation shows that

$\left( i\right) \ $The Runge - Kutta integrator with $1$- step is a Poisson
integrator if and only if%
\begin{equation*}
\left\{ 
\begin{array}{c}
h=\frac{1}{2bA-aA},\ \ \ \ \  \\ 
2bA-aA\neq 0.%
\end{array}%
\right.
\end{equation*}

$\left( ii\right) \ $The Runge - Kutta integrator with $1$- step is an
energy integrator if and only if%
\begin{equation*}
A=0\text{, sau }B=0\text{.}
\end{equation*}

$\left( iii\right) \ $The Gauss-Legendre integrator is a Poisson
integrator,

$\left( iv\right) \ $The Gauss-Legendre integrator is an energy
integrator.

\textbf{Proposition 3.1. }([8]) Let $\left( 
\mathbb{R}
^{n},\left\{ \cdot ,\cdot \right\} ,H\right) $ be a Hamilton-Poisson mechanical
system. If the matrix $\Pi $ is constant, then the Runge-Kutta integrator
with $s$-steps is a Poisson integrator if and only if%
\begin{equation*}
\left\{ 
\begin{array}{c}
b_{i}a_{ij}+b_{j}a_{ji}-b_{i}b_{j}=0, \\ 
1\leq i,\ j\leq s.\ \ \ \ \ \ \ \ \ \ \ \ \ \ \ \ \ 
\end{array}%
\right. .
\end{equation*}

\textbf{Observation 3.2. }In the particular case $n=2m\ \text{and }\Pi =%
\left[ 
\begin{array}{cc}
0_{m} & I_{m} \\ 
-I_{m} & 0_{m}%
\end{array}%
\right] \ $the theorem [7], [15] and [17] is found (see \textbf{%
Proposition 2.2}).

\textbf{Proposition 3.2. }([1]) Let $\left( 
\mathbb{R}
^{n},\left\{ \cdot ,\cdot \right\} ,H\right) $ be a Hamilton-Poisson
mechanical system. If the matrix $\Pi $ is constant, then the Gauss-Legendre
integrator is a Poisson integrator.

\textbf{Observation 3.3.} In the particular case $n=2m\ \text{and }\Pi =%
\left[ 
\begin{array}{cc}
0_{m} & I_{m} \\ 
-I_{m} & 0_{m}%
\end{array}%
\right] \ $the theorem [5] is found (see \textbf{%
Proposition 2.1.}).

\bigskip

\textbf{4. The Lie-Trotter Integrator and some properties in the dynamics
of symmetric free rigid body }

\ \ \ 

Let $\left( 
\mathbb{R}
^{n},\left\{ \cdot ,\cdot \right\} ,H\right) $ be a Hamilton-Poisson
mechanical system. The Lie-Trotter integrator (see [18]) usually applies when the Hamiltonian $H$ can be written as%
\begin{equation*}
H=H_{1}+H_{2},
\end{equation*}%
and the dynamics generated by $X_{H_{1}}$ $\left[ \text{resp }X_{H_{2}}%
\right] $, hence $\exp \left( tX_{H_{1}}\right) $ $\left[ \text{resp. }\exp
\left( tX_{H_{2}}\right) \right] $ can be explicitly integrated, hence $\exp
\left( tX_{H_{1}}\right) $ and $\exp \left( tX_{H_{2}}\right) $ can be
explicitly calculated. Then the Lie-Trotter integrator is given by%
\begin{equation}
\phi _{t}\left( x\right) =\left[ \exp \left( tX_{H_{2}}\right) \circ \exp
\left( tX_{H_{1}}\right) \right] \left( x\right) .  \tag{$4.1$}
\end{equation}

Because $\phi _{t}$ is actually composing two Hamiltonian flows, we have

\textbf{Proposition 4.1. }([12], [13]) The Lie-Trotter integrator $\left( 4.1\right) $
has the following properties

$\left( i\right) \ \left\{ \phi _{t}\right\} _{t\in 
\mathbb{R}
}$ is a Poisson integrator,

$\left( ii\right) \ $The restriction of the foliation of a manifold
symplectic Poisson $\left( 
\mathbb{R}
^{n},\left\{ \cdot ,\cdot \right\} \right) $ gives rise to a symplectic
integrator.

\bigskip

Obviously the integrator $(4.1)$ is a first order integrator. This order can be
grown by using some ideas from [2]. Concretely, a Poisson integrator of second order is given as%
\begin{equation*}
\phi _{t}^{2}\left( x\right) =\left[ \exp \left( \frac{t}{2}X_{H_{1}}\right)
\circ \exp \left( tX_{H_{2}}\right) \circ \exp \left( \frac{t}{2}%
X_{H_{1}}\right) \right] \left( x\right) .
\end{equation*}

The 4th order Poisson integrator is given by%
\begin{equation*}
\phi _{t}^{4}\left( x\right) =\left[ \phi _{t}^{2}\left( x_{1}t\right) \circ
\phi _{t}^{2}\left( x_{0}t\right) \circ \phi _{t}^{2}\left( x_{1}t\right) %
\right] \left( x\right) ,
\end{equation*}%
where%
\begin{equation*}
x_{0}=\frac{\sqrt[3]{2}}{\sqrt[3]{2}-2},\ x_{1}=\frac{1}{2-\sqrt[3]{2}}.
\end{equation*}

The process continues and, finally, it is obtained a Poisson integrator of the
even order $2n+2$ defined as%
\begin{equation*}
\phi _{t}^{2n+2}\left( x\right) =\left[ \phi _{t}^{2n}\left( x_{1}t\right)
\circ \phi _{t}^{2n}\left( x_{0}t\right) \circ \phi _{t}^{2n}\left(
x_{1}t\right) \right] \left( x\right) ,
\end{equation*}%
where%
\begin{equation*}
x_{0}=\frac{\sqrt[2n+1]{2}}{\sqrt[2n+1]{2}-2},\ x_{1}=\frac{1}{2-\sqrt[2n+1]{%
2}}.
\end{equation*}

Analogously explicit Poisson integrators can be build unconcerned their even
order. We note that construction is difficult because in the construction of
a Poisson integrator of even order $2n$ , the second-order integrator $\phi
_{t}^{2}$ is used for $3^{n-1}$ times. Therefore the number of steps $k$ is $%
1+3^{n-1}$, hence%
\begin{equation*}
k=1+3^{n-1},
\end{equation*}%
which grows extremely fast, and the calculations become increasingly more
complex.

We proof below the next theorem:

\textbf{Theorem 4.1. }If the Lie - Trotter integrator is obtained, in the
case of Euler equations for the dynamics of symmetric free rigid body, then
it is a Poisson integrator.

\textbf{Proof}

The Euler equations for the dynamics of symmetric free rigid body are
written as%
\begin{equation}
\left\{ 
\begin{array}{l}
\overset{\cdot }{m}_{1}=a_{1}m_{2}m_{3,} \\ 
\overset{\cdot }{m}_{2}=-a_{1}m_{1}m_{3}, \\ 
\overset{\cdot }{m}_{3}=0,%
\end{array}%
\right.  \tag{4.2}
\end{equation}%
where $a_{1}=\dfrac{1}{I_{3}}-\dfrac{1}{I_{1}}$, the $I_{1}$, $I_{3}$
inertia tensors are components of the body and we assume in all that follows
that $I_{1}>I_{3}>0$.

As we well know (see some Puta's papers) the dynamic $\left( 4.2\right) $ has the following Hamilton-Poisson realization $\left(
\left( so\left( 3\right) \right) ^{\ast }\simeq 
\mathbb{R}
^{3},\Pi ,H\right) $, where
\begin{equation*}
\Pi =\left[ 
\begin{array}{ccc}
0 & -m_{3} & m_{2} \\ 
m_{3} & 0 & -m_{1} \\ 
-m_{2} & m_{1} & 0%
\end{array}%
\right] ,\ H\left( m_{1},m_{2},m_{3}\right) =\frac{1}{2}\left( \frac{%
m_{1}^{2}}{I_{1}}+\frac{m_{2}^{2}}{I_{1}}+\frac{m_{3}^{2}}{I_{3}}\right) .
\end{equation*}

We can say more and more precisely that the function $C\in C^{\infty }\left( 
\mathbb{R}
^{3},%
\mathbb{R}
\right) $ given by%
\begin{equation*}
C\left( m_{1},m_{2},m_{3}\right) =\frac{1}{2}\left(
m_{1}^{2}+m_{2}^{2}+m_{3}^{2}\right) ,
\end{equation*}%
is one Casimir of our configuration $\left( 
\mathbb{R}
^{3},\Pi \right) .$ Actually, the symplectic foliation of the Poisson variety $%
\left( 
\mathbb{R}
^{3},\Pi \right) $ is given by $\left\{ 0\right\} $ and by the coadjuncte orbits $%
\left( O_{k},\omega _{k}\right) ,$ where%
\begin{equation*}
O_{k}=\left\{ \left( m_{1},m_{2},m_{3}\right) \in 
\mathbb{R}
^{3}\mid m_{1}^{2}+m_{2}^{2}+m_{3}^{2}=k^{2}\right\} ,
\end{equation*}%
and%
\begin{equation*}
\omega _{k}=\frac{1}{k}\left( m_{2}dm_{1}\wedge dm_{3}-m_{3}dm_{1}\wedge
dm_{2}-m_{1}dm_{2}\wedge dm_{3}\right) ,
\end{equation*}%
$\left( \omega _{k}\text{ is symplectic form Kirillov-Kostant-Souriau}%
\right) .$

It can be now seen that the Hamiltonian field $X_{H}$ splits%
\begin{equation*}
X_{H}=X_{H_{1}}+X_{H_{2}}+X_{H_{3}},
\end{equation*}%
where%
\begin{equation*}
H_{1}\left( m_{1},m_{2},m_{3}\right) =\frac{1}{2I_{1}}m_{1}^{2},\
H_{2}\left( m_{1},m_{2},m_{3}\right) =\frac{1}{2I_{1}}m_{2}^{2},\
H_{3}\left( m_{1},m_{2},m_{3}\right) =\frac{1}{2I_{3}}m_{3}^{2}.
\end{equation*}

Then the corresponding flows have the following expressions 
\begin{equation*}
\left[ 
\begin{array}{c}
m_{1}\left( t\right) \\ 
m_{2}\left( t\right) \\ 
m_{3}\left( t\right)%
\end{array}%
\right] =\left[ 
\begin{array}{ccc}
1 & 0 & 0 \\ 
0 & \cos \frac{m_{1}\left( 0\right) }{I_{1}}t & \sin \frac{m_{1}\left(
0\right) }{I_{1}}t \\ 
0 & -\sin \frac{m_{1}\left( 0\right) }{I_{1}}t & \cos \frac{m_{1}\left(
0\right) }{I_{1}}t%
\end{array}%
\right] \left[ 
\begin{array}{c}
m_{1}\left( 0\right) \\ 
m_{2}\left( 0\right) \\ 
m_{3}\left( 0\right)%
\end{array}%
\right] ,
\end{equation*}%
\begin{equation*}
\left[ 
\begin{array}{c}
m_{1}\left( t\right) \\ 
m_{2}\left( t\right) \\ 
m_{3}\left( t\right)%
\end{array}%
\right] =\left[ 
\begin{array}{ccc}
\cos \frac{m_{2}\left( 0\right) }{I_{1}}t & 0 & -\sin \frac{m_{2}\left(
0\right) }{I_{1}}t \\ 
0 & 1 & 0 \\ 
\sin \frac{m_{2}\left( 0\right) }{I_{1}}t & 0 & \cos \frac{m_{2}\left(
0\right) }{I_{1}}t%
\end{array}%
\right] \left[ 
\begin{array}{c}
m_{1}\left( 0\right) \\ 
m_{2}\left( 0\right) \\ 
m_{3}\left( 0\right)%
\end{array}%
\right] ,
\end{equation*}%
\begin{equation*}
\left[ 
\begin{array}{c}
m_{1}\left( t\right) \\ 
m_{2}\left( t\right) \\ 
m_{3}\left( t\right)%
\end{array}%
\right] =\left[ 
\begin{array}{ccc}
\cos \frac{m_{3}\left( 0\right) }{I_{3}}t & \sin \frac{m_{3}\left( 0\right) 
}{I_{3}}t & 0 \\ 
-\sin \frac{m_{3}\left( 0\right) }{I_{3}}t & \cos \frac{m_{3}\left( 0\right) 
}{I_{3}}t & 0 \\ 
0 & 0 & 1%
\end{array}%
\right] \left[ 
\begin{array}{c}
m_{1}\left( 0\right) \\ 
m_{2}\left( 0\right) \\ 
m_{3}\left( 0\right)%
\end{array}%
\right] .
\end{equation*}

Then results that Lie - Trotter integrator has the following expression%
\begin{equation}
\left[ 
\begin{array}{c}
m_{1}^{n+1} \\ 
m_{2}^{n+1} \\ 
m_{3}^{n+1}%
\end{array}%
\right] =M\cdot N\cdot P\cdot \left[ 
\begin{array}{c}
m_{1}^{n} \\ 
m_{2}^{n} \\ 
m_{3}^{n}%
\end{array}%
\right] ,  \tag{4.3}
\end{equation}%
where%
\begin{equation*}
M=\left[ 
\begin{array}{ccc}
1 & 0 & 0 \\ 
0 & \cos \frac{m_{1}\left( 0\right) }{I_{1}}t & \sin \frac{m_{1}\left(
0\right) }{I_{1}}t \\ 
0 & -\sin \frac{m_{1}\left( 0\right) }{I_{1}}t & \cos \frac{m_{1}\left(
0\right) }{I_{1}}t%
\end{array}%
\right] ,
\end{equation*}%
\begin{equation*}
N=\left[ 
\begin{array}{ccc}
\cos \frac{m_{2}\left( 0\right) }{I_{1}}t & 0 & -\sin \frac{m_{2}\left(
0\right) }{I_{1}}t \\ 
0 & 1 & 0 \\ 
\sin \frac{m_{2}\left( 0\right) }{I_{1}}t & 0 & \cos \frac{m_{2}\left(
0\right) }{I_{1}}t%
\end{array}%
\right] ,
\end{equation*}%
\begin{equation*}
P=\left[ 
\begin{array}{ccc}
\cos \frac{m_{3}\left( 0\right) }{I_{3}}t & \sin \frac{m_{3}\left( 0\right) 
}{I_{3}}t & 0 \\ 
-\sin \frac{m_{3}\left( 0\right) }{I_{3}}t & \cos \frac{m_{3}\left( 0\right) 
}{I_{3}}t & 0 \\ 
0 & 0 & 1%
\end{array}%
\right] .
\end{equation*}

In conclusion the Lie-Trotter integrator $\left( 4.3\right) $ is a Poisson integrator.

From the above, (see [13]), we get

\textbf{Proposition 4.2}: The restriction of the coadjuncte orbits $\left(
O_{k},\omega _{k}\right) $ involves obtaining a symplectic integrator.

\textbf{Proof:} It results from the above theorem and \textbf{Proposition 4.1%
}.

\ \ \

\textbf{Remark 4.1}: The Lie-Trotter integrator is not an energy integrator.

\ \ \

\textbf{Remark 4.2}: If we compare with the Runge-Kutta integrator with $4$
steps, then we obtain almost the same results, but the Lie - Trotter
integrator has the advantage that it is much easier to implement.

\ \ \

\textbf{Remark 4.3} Similar results with the ones presented in this paper are shown in [13] for the free rigid body. By analysing [13], one can see that some of the computations are identical.

\ \ \

\textbf{Remark 4.4}: It remains an open problem to find if there are characteristic roots
of the matrix $M\cdot N\cdot P$ which have even values and follow a Hopf
bifurcation moving from the unit circle on the real right so that their product
is allways $1$ and if there is a transition in the dynamic integration $(4.2)$ by
means of the Lie - Trotter integrator, from the \textit{real world} in the
\textit{complex world} tagged with the Hopf bifurcation for some of
their values coming from $M\cdot N\cdot P$.

\ \ \

\textbf{Remark 4.5} The square matrix $\ R=M\cdot N\cdot P$ , where $\Pi
_{n}=R^{n}\cdot \Pi _{0}$, $\Pi _{0}=\left[ 
\begin{array}{c}
1 \\ 
0 \\ 
0%
\end{array}%
\right] $ can also give an open problem for finding if there are values of $n
$ so that the Markov chain related to the matrix can be found in a maximum
possible number of states.

\bigskip
\bigskip
\bigskip

\bigskip
\bigskip
\bigskip
\textbf{Acknowledgement}

This year we mark a decade from the death of the regretted Professor Dr. Mircea Puta, the coordinator of my PhD thesis. I dedicate this paper in his memory and for the strong collaboration with him.  
\bigskip

\textbf{References}

[1]  Austin, M. A., Krishnaprassad P.S., Wang L.S., Almost Poisson integration
of rigid body systems, Journal of Computational Physics, vol. 107, No. 1 (1993),
105-117.

[2] Benzel, S., Ge, Z., Scovel, C., Elementary construction of higher order
Lie-Poisson integrators, Phys. Letter A, 174 (1993), 229-232.

[3] Buchner, K., Craioveanu, M. and Puta, M., Recent progress in the integration
of Poisson systems via de mid-point rule and Runge-Kutta algorithm, Balkan
Journal of Geometry and Its Applications, Vol. 1. Nr. 20 (1996), 9-20.

[4] Channel, P. and Scovel, C., Symplectic integration of Hamiltonian systems,
Nonlinearity 3 (1990), 231-259.

[5] Feng, K., Difference schemes for Hamiltonian formalism and symplectic
geometry, J. Comp. Math. 4 (1985), 279-289.

[6] Ge, Z. and Marsden, J., Lie-Poisson integrators and Lie-Poisson
Hamiltonian-Jacobi theory, Phys. Lett. A 133 (1988), 134-139.

[7] Lasagni, F., Canonical Runge-Kutta methods, ZAMP 39 (1988), 952-953.

[8] Mc Lachlan, R.L., Comment on Poisson schemes for Hamiltonian systems on
Poisson manifolds, Computers Math. Applic. vol. 29 No. 3 $\left(
1995\right) .$

[9] Mo\c{s}, I., An Introduction to Geometric Mechanics, Cluj University
Press. $\left( 2005\right) .$

[10] Neri, F.,\ Lie-algebras and canonical integration, preprint, Dept. of
Physics, University of Maryland, $\left( 1987\right) .$

[11] Puta, M., Poisson integrators, Analele Univ. Timi\c{s}oara, vol. 31 (1993),
267-272.

[12] Puta, M., An overview of some Poisson integrators, in Enumath. 97, 2nd
European Conference on Numerical Mathematics and Advanced Applications, eds.
H.G. Bock, G. Kanschat, Y. Kuznetsov and J. Periaux, World Scientific,(1998),
518-523.

[13] Puta, M., Lie-Trotter formula and Poisson dynamics, Int. Journ. of
Biffurcation and Chaos, vol. 9 (1999), 555-559.

[14] Ruth, R.D., A canonical integration technique, IEEE Trans. Nucl. Sci.,
N.S. 30 (1983), 2669-2271.

[15] Sanz-Serna, J.M., Runge-Kutta schemes for Hamiltonian systems, BIT 28 (1988),
877-883.

[16] Sanz-Serna, J.M., Symplectic integrators for Hamiltonian problems an
overview, Acta Numerica 1,(1992), 243-286.

[17] Suris, Y.B., Canonical transformations generated by methods of Runge-Kutta
type for the numerical integration of the system $\overset{\cdot \cdot }{x}%
=\partial U/\partial x$ (Russian), Zh. Vychisl. Mat. i Mat. Fiz. 29 (1989),
202-211.

[18] Trotter, H.F., On the product of semigroups of operators, Proc. Amer.
Math. Soc. 10 (1959), 545-551.

[19] Yoshida, H., Construction of higher order symplectic integrators, Phys.
Lett. A, 150 (1990), 262-268.

\bigskip

\end{document}